\documentclass[12pt]{article}
\usepackage{amsthm,amsmath,amssymb,latexsym,graphics}
\newtheorem{proposition}{Proposition}
\newtheorem{theorem}{Theorem}
\newtheorem{remark}{Remark}

\newcommand{\grad}{\nabla}

\begin{document}
\title{Classification of singular radial solutions  to the $\sigma_k$
Yamabe equation on annular domains}
\author{S.-Y. Alice  Chang \thanks{partially supported  by  NSF through grant DMS0245266.} \\
Department of Mathematics \\Fine Hall\\ Princeton University \\Princeton, NJ 08540\\ \texttt{chang@math.princeton.edu}
\and
Zheng-Chao Han \thanks{corresponding author; partially supported  by  NSF through grant DMS-0103888.}
\\Department of Mathematics \\Rutgers University \\110 Frelinghuysen Road \\
Piscataway, NJ 08854\\ \texttt{zchan@math.rutgers.edu} 
\and
Paul Yang \thanks{partially supported  by  NSF through grant DMS0245266.}\\
Department of Mathematics \\Fine Hall \\Princeton University \\Princeton, NJ 08540\\ 
\texttt{yang@math.princeton.edu}}
\date{}
\maketitle
\begin{abstract}
The study of the $k$-th elementary symmetric function of the Weyl-Schouten curvature tensor of a Riemannian metric,
the so called $\sigma_k$ curvature,
has produced many fruitful results in conformal geometry in recent years. In these studies in conformal
geometry, the deforming conformal factor is considered to be a solution of  a fully nonlinear elliptic PDE.
Important  advances have been made in recent years in the understanding of the analytic behavior of solutions
of the PDE. However, the singular behavior of these solutions, which is important in describing many
important questions in conformal geometry, is little understood. This note classifies
 all possible radial solutions, in particular, the
\emph{singular}  solutions of the $\sigma_k$ Yamabe equation, which describes conformal metrics whose 
$\sigma_k$ curvature equals a constant. Although the analysis involved is of elementary nature,
these results should provide useful guidance
in studying the behavior of singular solutions in the general situation. 
\end{abstract}


\section{Description of the results}
This note is concerned with radial solutions to the equation 
\begin{equation} \label{1}
\sigma_k (A_g) = \text{constant},
\end{equation}
on domains of the form $\{\,x\in \mathbb R^n : r_1 < |x| < r_2 \,\}$,
 where $A_g$ is the
Weyl-Schouten tensor of the conformal metric $g = v^{-2}(|x|) |dx|^2$: 
\[
A_g =  \frac {1}{n-2} \{ Ric - \frac {R}{2(n-1)}g \},
\]
$\sigma_k(A_g)$ denotes the $k$-th elementary symmetric function of the eigenvalues of $A_g$
with respect to $g$, and $0 \le r_1 < r_2 \le \infty$. 
In this note we classify all possible radial solutions, in particular, the
\emph{singular} solutions. The main motivation for our study is to use these results   as guidance
in studying the behavior of singular solutions in the general situation. 

There are at least two kinds of 
geometric considerations that lead to the study of singular solutions
of equations of the above type.
The first kind is related to the characterization of the size of 
the limit set of the image domain in $\mathbb S^n$ of
the developing map of a locally conformally flat $n$-manifold. More specifically, one is led to  find  
necessary/sufficient conditions  on a domain $\Omega \subset
\mathbb S^n$ so that it 
admits a conformal metric $g = v^{-2}(x) |dx|^2$ which is complete, 
and with its Weyl-Schouten
tensor $A_g$ in the $\Gamma^{\pm}_k$ class, \emph{i.e.}, the eigenvalues, 
$\lambda_1 \le \cdots \le \lambda_n$, of $A_g$ at each $x \in  \Omega$ satisfy
$ \sigma_j (\lambda_1 , \cdots , \lambda_n) >0$ for all $j, 1\le j \le k$, 
in the case of 
$\Gamma^+_k$; and $(-1)^j \sigma_j (\lambda_1 , \cdots , \lambda_n) >0$ 
for all $j, 1\le j \le k$, in the case of $\Gamma^-_{k}$.  For $k \ge
2$, it is natural to restrict to metrics whose Weyl-Schouten tensor is in  the $\Gamma^{\pm}_k$ class,
 because, for a metric in such a class, \eqref{1} becomes a  \emph{fully nonlinear} PDE in $v$ that is 
\emph{elliptic}~.
 In the case of $k=1$, $\sigma_1(A_g)$ is simply a constant multiple of the scalar curvature of $g$; so
 $A_g$ in the $\Gamma^{\pm}_k$ class is a generalization of the notion that the scalar
curvature $R_g$ of $g$ having a fixed $\pm$ sign. For the positive scalar curvature case, 
in \cite{SY}, Schoen-Yau proved 
 that if a complete metric $g = v^{-2}(x) |dx|^2$
exists on a domain $\Omega \subset \mathbb S^n$ with $\sigma_1 (A_g)$ having a positive lower bound, 
then the Hausdorff dimension of $\partial \Omega$ has to be $\leq \frac{n-2}{2}$. Later 
Mazzeo and Parcard \cite{MP1} proved that if $\Omega \subset \mathbb S^n$ is a domain such that
$\mathbb S^n \setminus \Omega$ consists a finite number of smooth submanifolds of dimension
 $1\le k \le \frac{n-2}{2}$,
then one can find a complete metric $g = v^{-2}(x) |dx|^2$ on  $ \Omega$ with its scalar
curvature idential to $+1$. For the negative scalar curvature case, the work of
 Lowner-Nirenberg \cite{LN}, Aviles \cite{Av}, and Veron \cite{Ver} implies that if $\Omega \subset \mathbb S^n$
admits a complete, conformal metric with negative constant scalar curvature, then the Hausdorff
dimension of $\partial \Omega > \frac{n-2}{2}$.   Lowner-Nirenberg \cite{LN} also proved that  if 
$\Omega \subset \mathbb S^n$ is a domain with smooth boundary $\partial \Omega$ of dimension $> \frac{n-2}{2}$,
then there exists a  complete metric $g = v^{-2}(x) |dx|^2$
on  $\Omega$ with $\sigma_1 (A_g) = -1$. This result was later generalized by D. Finn \cite{Fin} to the case
of $\partial \Omega$ consisting of smooth submanifolds of dimension $> \frac{n-2}{2}$ and  with boundary. 
So the consideration of singular solutions of equations of type \eqref{1} can be considered as a natural 
generalization of these known results. In fact, in \cite{CHgY}, Chang, Hang, and Yang proved that if
$\Omega \subset \mathbb S^n$ ($n \geq 5$) admits a complete, conformal metric $g$ with 
$$\sigma_1(A_g) \geq c_1 >0, \quad  \sigma_2(A_g) \geq 0, \quad \text{and}$$
\begin{equation} \label{ub}
|R_g| +|\grad_g R|_g \leq c_0, 
\end{equation}
then $\dim (\mathbb S^n \setminus \Omega) < \frac{n-4}{2}$. This has been generalized by M. Gonzalez to the
case of $2 < k < n/2$: if
$\Omega \subset \mathbb S^n$  admits a complete, conformal metric $g$ with 
$$\sigma_1(A_g) \geq c_1 >0, \quad  \sigma_2(A_g),\;\cdots\; \sigma_k(A_g)  \geq 0, \quad \text{and}$$
\eqref{ub}, then $\dim (\mathbb S^n \setminus \Omega) < \frac{n-2k}{2}$. See also the work of Guan, Lin and 
Wang \cite{GLW}.

From a more broad perspective, one main reason for the attention to the $\sigma_k$-curvature ($k >1$)
in conformal geometry is that these curvatures place a much stronger control on the curvature tensor.
In dimension $4$, the $\sigma_2$ curvature comes into play in the Chern-Gauss-Bonnet formula. Chang, 
Gursky and Yang 
\cite{CGY1} observed that if $\sigma_1(A_g), \sigma_2(A_g) >0$ at a point on a $4$-dimensional manifold, then
the Ricci tensor of $g$ is positive definite at that point. This algebraic relation has been generalized to
higher dimensions by Guan, Viaclovsky, and Wang \cite{GVW}. The first important application of the
$\sigma_k$-curvature to conformal geometry is the main theorem in \cite{CGY1}, where the authors proved that
if (i) $\int_{M^4} \sigma_2 (A_g) d\;vol$, which is conformally invariant on $M^4$, is  positive; and (ii)
the Yamabe class of $(M^4, g)$ is positive, then there is a conformal metric $\tilde g = e^{2w}g$ on 
$M^4$ such that $A_{\tilde g} \in \Gamma_2^+$.

The second kind of consideration for studying singular 
solutions of \eqref{1} 
is due to a basic phenomenon of solutions of the more general equation, 
$ \sigma_k (A_g) = f(x, v)$:
for a solution $v$ with its  $A_g \in \Gamma^+_k$, there is  higher
derivative estimates for the solution $v$ in
terms of the $C^0$ estimates of $v$ \cite{Via1} \cite{GW1}; while for a solution $v$
with its $A_g \in \Gamma^-_k$, there is a lack of second derivative
estimates for $v$, even when the $C^0$ norm of $v$ is under control. 
The singular radial solutions exhibit this behavior explicitly.

Let us formulate our tasks more explicitly:
the problem of the classification of radial solutions of \eqref{1}
 consists of determining (a) the maximal domain of definition of
 each solution (a finite ball, a finite punctured ball,
an annulus,  $\mathbb R^n \setminus \{0\}$,  or the entire $\mathbb R^n$?); (b) the limiting
behavior of each  solution upon approaching the limits of its domain of
definition and the geometric meanings of such behavior
 (completeness vs incompleteness, etc); 
and (c) whether the solution is the the $\Gamma^+$ or $\Gamma^-$ class.

First, let us work out \eqref{1} more explicitly in the case of radial 
solutions. Let 
$$g = v^{-2}(|x|) |dx|^2, \quad \text{and} \quad |x| = r.$$
Then 
\begin{eqnarray*}
A_{ij}  &= \frac {v_{ij}}{v} - \frac {|\grad v|^2}{2v^2} \delta_{ij} \\
        &= \lambda \delta_{ij} + \mu \frac{x_ix_j}{|x|^2},
\end{eqnarray*}
with $\lambda = \frac{v_r}{rv} (1 - \frac{rv_r}{2v})$ and 
 $\mu = \frac{v_{rr}}{v} -\frac{v_r}{rv}$.
The eigenvalues of $A$ with respect to $|dx|^2$ are $\lambda$ with
 multiplicity $(n-1)$, and 
$\lambda +\mu$ with multiplicity $1$. The formula for $\sigma_k(A_g)$ can
be found easily by the binomial expansion of $(x-\lambda)^{n-1}(x-\lambda-\mu)$:
\begin{equation} \label{2}
\sigma_k(A_g) = c_{n,k} v^{2k} \lambda^{k-1} (n\lambda + k \mu),
\end{equation}
where $c_{n,k} = \frac {(n-1)!}{k! (n-k)!}$.

First, the following   observation follows directly from \eqref{2}:

\begin{remark}
If $k> 1$ and $v(r)>0 $ is a $C^2$ function on $r_1 < r < r_2$
 such that  $\sigma_k(A_{g_v})$
 has  a fixed sign for $r_1 < r < r_2$, then $v(r)$ 
is strictly monotone on $r_1 < r < r_2$. More
specifically $\lambda= \frac{v_r}{rv} (1 - \frac{rv_r}{2v})$ never changes sign
 on $r_1 < r < r_2$ .
\end{remark}

Next we introduce new variables, $t = \ln r$,  $x = r \omega$ with
 $\omega \in \mathbb S^{n-1}$, and $\xi(t)$ such that
$$g= v^{-2}(r) |dx|^2 = e^{-2\xi} (dt^2 +d \omega^2),$$
then 
\[
\xi + t = \ln v , \qquad v_r=e^{\xi}(\xi_t+1)=(\xi_t+1)v/r, 
\]
and
\[
v_{rr}= e^{\xi-t}[\xi_{tt}+\xi_t(\xi_t+1)] = [\xi_{tt}+\xi_t(\xi_t+1)] ve^{-2t}.
\]
Thus
\[
\lambda = \frac {1-\xi_t^2}{2e^{2t}}, \qquad \text{and} \qquad 
\mu = e^{-2t}(\xi_{tt} + \xi_t^2 -1).
\]
 \eqref{2} then becomes
\begin{align} \label{3}
\sigma_k(A_g) &= c_{n,k} e^{2k(\xi+t)}\frac{(1-\xi_t^2)^{k-1}}{2^{k-1} e^{2(k-1)t}}\left[n \frac{1-\xi_t^2}{2e^{2t}} +  k\frac{\xi_{tt}+\xi_t^2-1}{e^{2t}}\right] \notag \\
&= c_{n,k}^{'} (1-\xi_t^2)^{k-1} \left[\frac{k}{n} \xi_{tt} + (\frac{1}{2} -\frac{k}{n})(1-\xi_t^2)\right]e^{2k\xi},
\end{align}
where $c_{n,k}^{'} = n c_{n,k} 2^{1-k}= 2^{1-k}\binom{n}{k} $.

The next proposition summarizes some further elementary properties of
solutions of \eqref{3}.
\begin{proposition} Consider any $C^2$ solution of \eqref{3} on an
interval $(t_{-}, t_{+})$  with $k>1$.
\begin{enumerate}
\item  If  $\sigma_k$ has a fixed
sign on the interval $(t_{-}, t_{+})$, then  either $1-\xi_t^2 > 0$ or
$1-\xi_t^2 < 0$ on $(t_{-}, t_{+})$. 
\item If $\sigma_k>0$ and $1 - \xi_t^2 >0$ on  $(t_{-}, t_{+})$, then
$v^{-2}(|x|)|dx|^2$ automatically stays in the $\Gamma_k^+$ class, 
\emph{i.e.}, it also satisfies $\sigma_l >0$ for any $1\leq l <k$.
\item If $\sigma_k>0$ and $1 - \xi_t^2 <0$ on  $(t_{-}, t_{+})$, and
 $k$ is even, then $v^{-2}(|x|)|dx|^2$ 
 automatically stays in the $\Gamma_k^-$ class, \emph{i.e.}, it also
 satisfies $(-1)^l \sigma_l>0$ for any $1 \leq l <k$.
\item   If  $\sigma_k<0$ and $1 - \xi_t^2 <0$ on  $(t_{-}, t_{+})$, and
$k$ is odd, then $v^{-2}(|x|)|dx|^2$
 automatically stays in the $\Gamma_k^-$ class.
\item  if $k>1$, $0< r_* < \infty$ is a limit point of the domain of definition of
a solution $v$ of \eqref{2} with $\sigma_k(A_g)=$ constant, and $v$ stays bounded away from $0$ and $\infty$
upon approaching $r_*$, then $v$ approaches a positive finite limit, and
its first derivative has either $v_r \to 0$, or $rv_r/v \to 2$,
 but its second derivative $v_{rr}$ blows up.
\end{enumerate}
\end{proposition}

Equation \eqref{3}, for $\sigma_k \equiv \text{constant}$, has a first
integral, \emph{i.e.}, a conserved quantity which reduces \eqref{3} to a first order ODE. 
After we completed our work, we realized that
in the case of $2k \neq n$, this fact  was already
pointed out by Viaclovsky in \cite{Via1} based on his variational
characterization of the solutions.  In fact, the first integral for
\eqref{3}, when $\sigma_k$ is a constant, 
can be obtained in all cases in a straightforward manner: simply 
multiplying both sides of \eqref{3} by $2ne^{-n\xi}\xi_t$, one
has, assuming $\sigma_k$ is normalized to be $2^{-k}\binom{n}{k}$, 
\[
[e^{(2k-n)\xi}(1-\xi_t^2)^k - e^{-n\xi}]_t \equiv 0.
\]
So $e^{(2k-n)\xi}(1-\xi_t^2)^k - e^{-n\xi}$ is a constant along any
solution. The remaining analysis of the bebavior of the radial solutions is
elementary: one uses the first integral to investigate the maximum
domain of definition of each solution and its asymptotic behavior upon
approaching the end points of its domain; one could also use the phase
plane portrait of \eqref{3} in the $\xi$-$\xi_t$ plane as a guidance for the global behavior of 
solutions of \eqref{3}. The phase plane portrait of \eqref{3}, when $\sigma_k$ is a constant,  depends on the
relation between $2k$ and $n$, as well as on whether $k$ is odd or even. 

 To help understand the different
situations that can possibly occur, the patterns of the phase portraits of \eqref{3}, when 
$\sigma_k$ is a \emph{positive} constant, are displayed on the next page. Note that the $x$ and $y$ 
axes in the displayed phase portraits stand for  the $\xi$ and $\xi_t$ axes, respectively.
\begin{figure}[hp!]
\scalebox{0.40}{\includegraphics{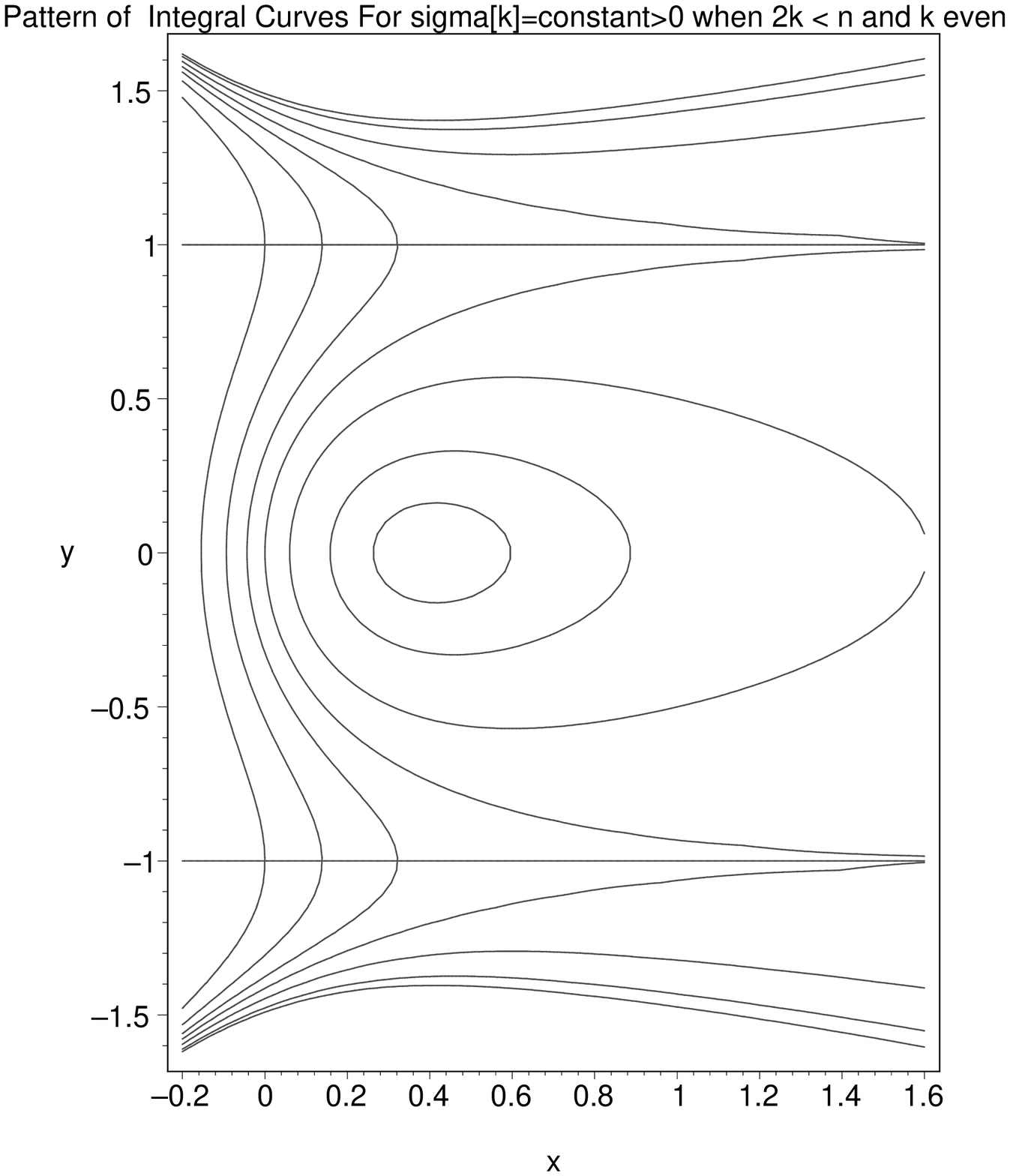}}\label{2ksnevenk}
\scalebox{0.40}{\includegraphics{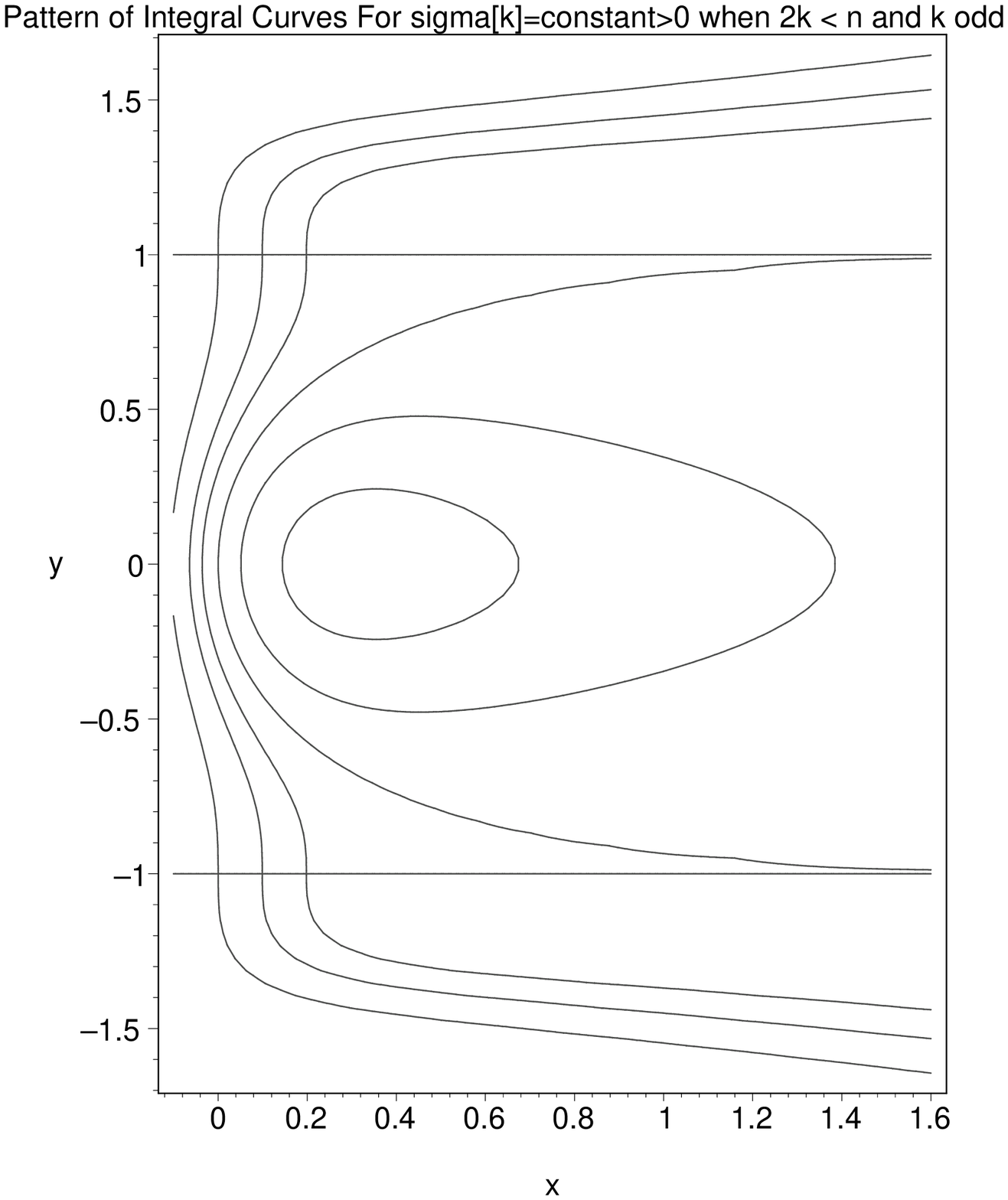}}\label{2ksnoddk}
\scalebox{0.40}{\includegraphics{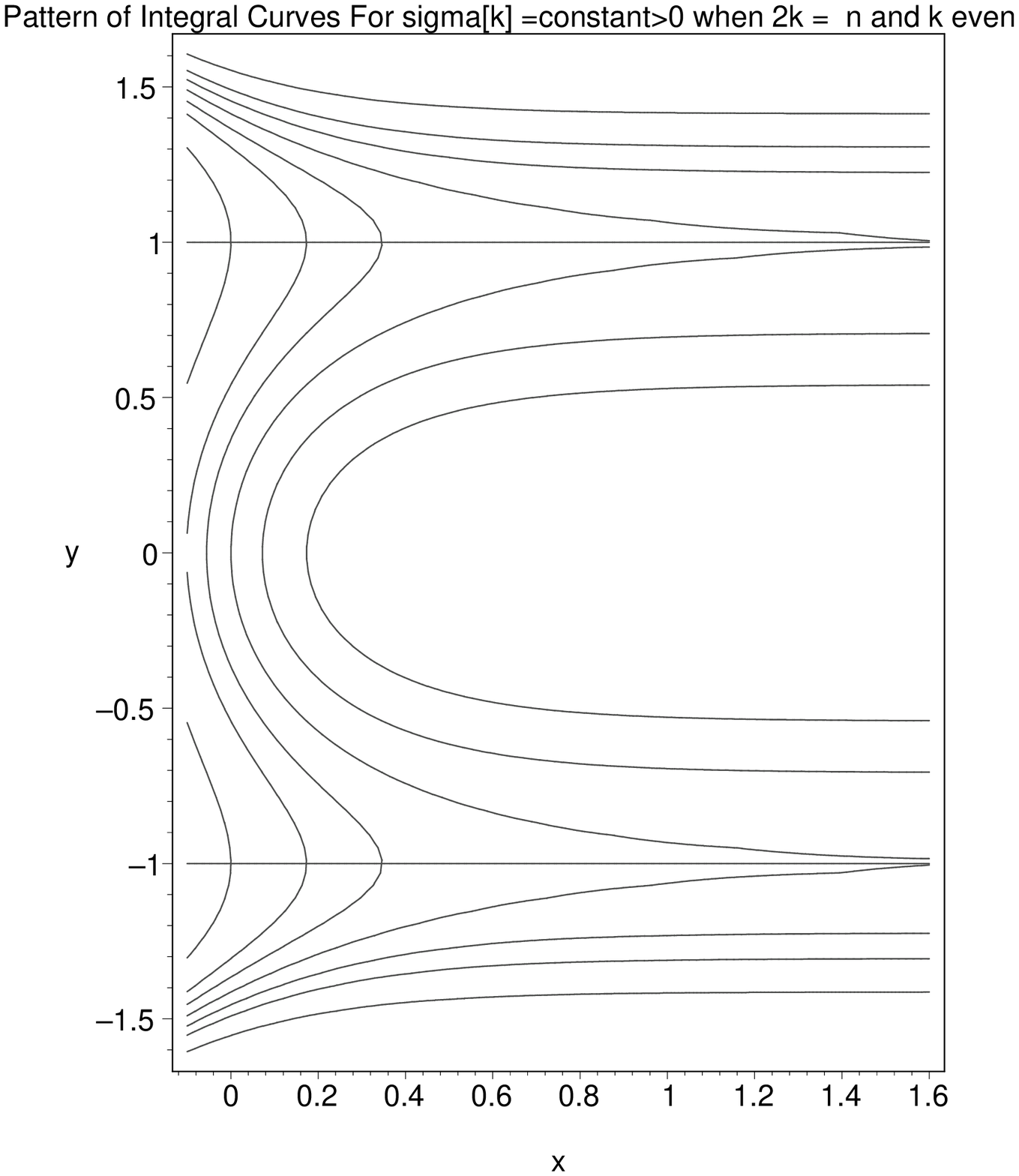}}\label{2kenevenk}
\scalebox{0.40}{\includegraphics{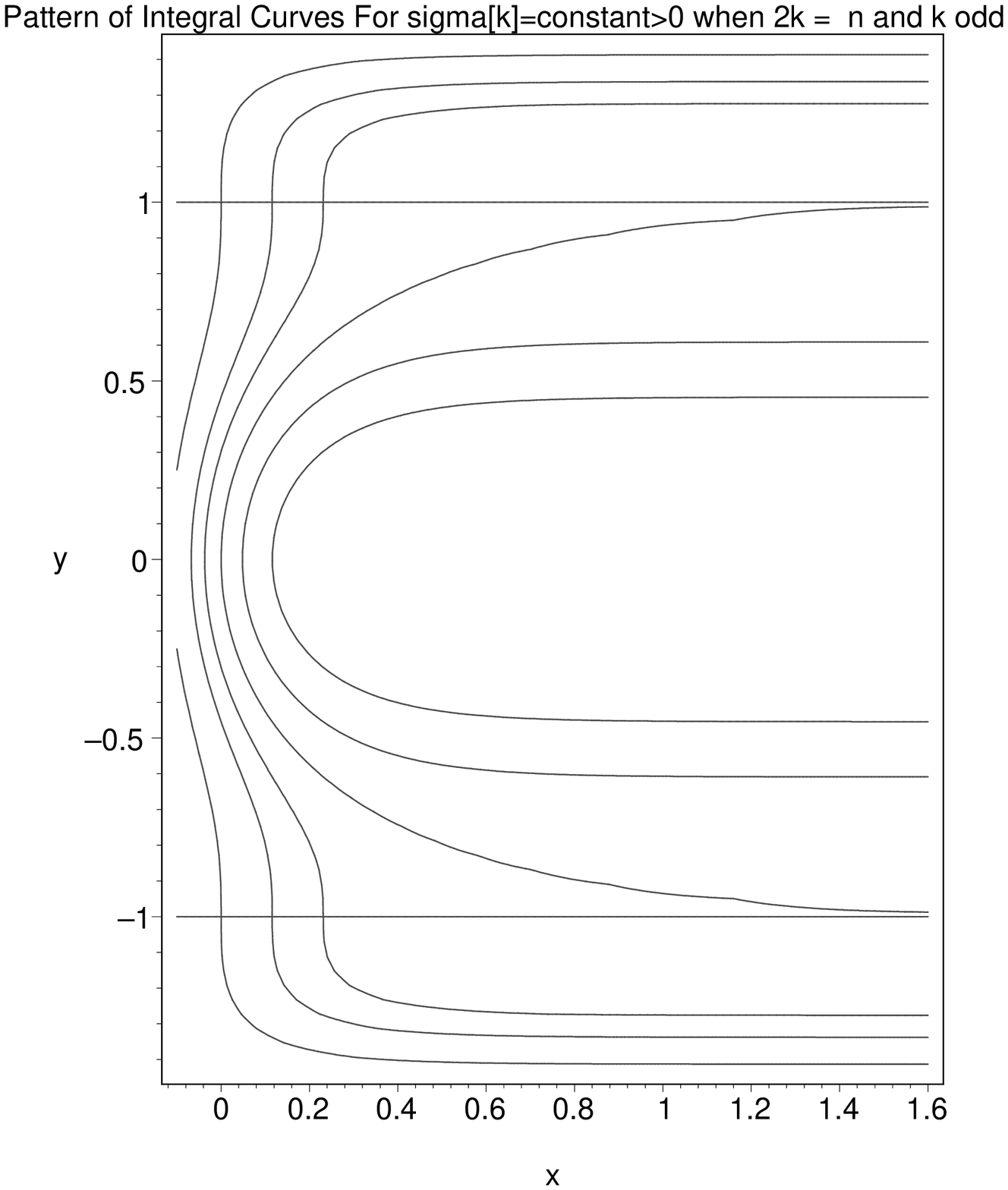}}\label{2kenoddk}
\scalebox{0.40}{\includegraphics{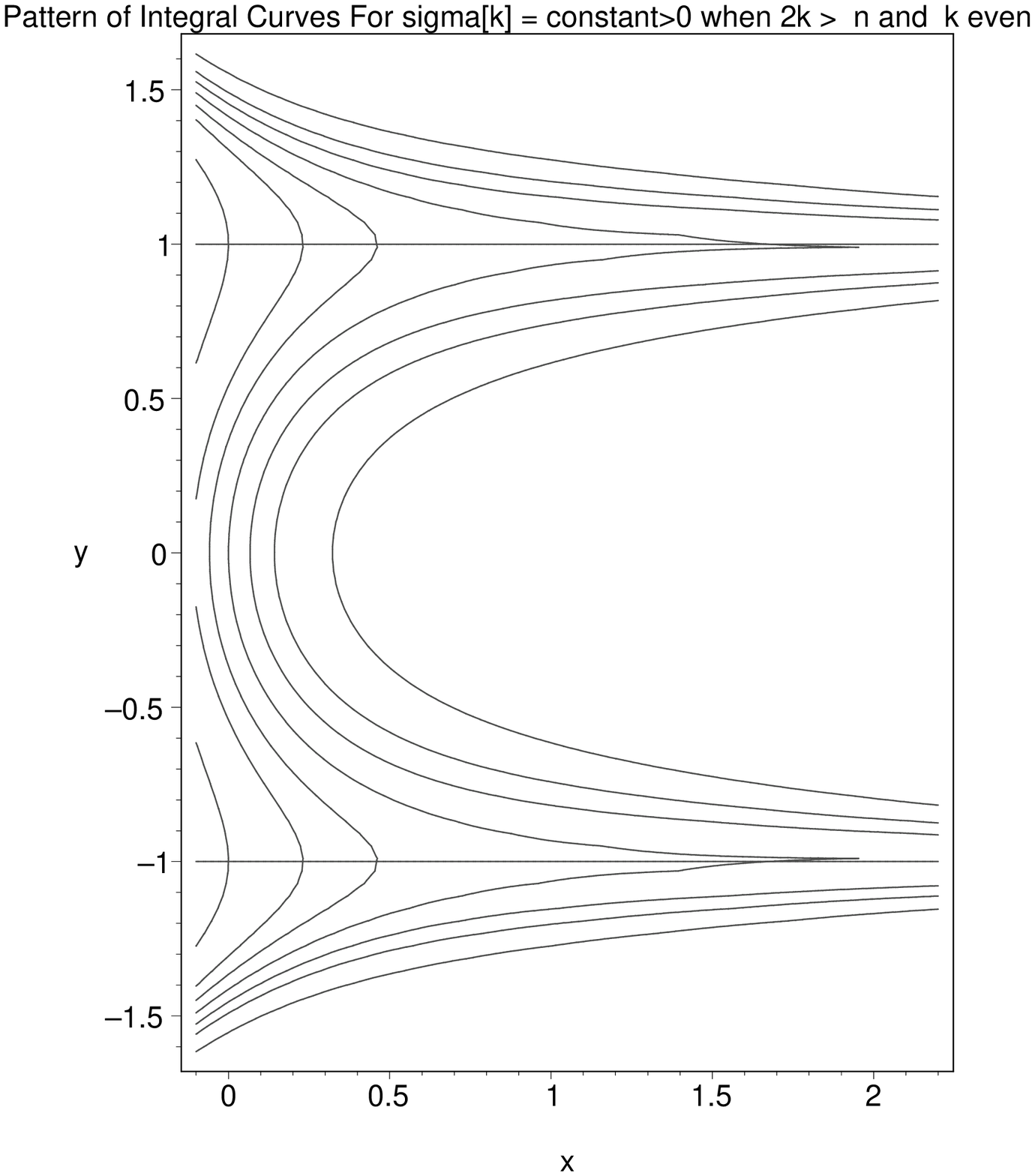}}\label{2kgnevenk}
\quad \qquad \qquad
\scalebox{0.40}{\includegraphics{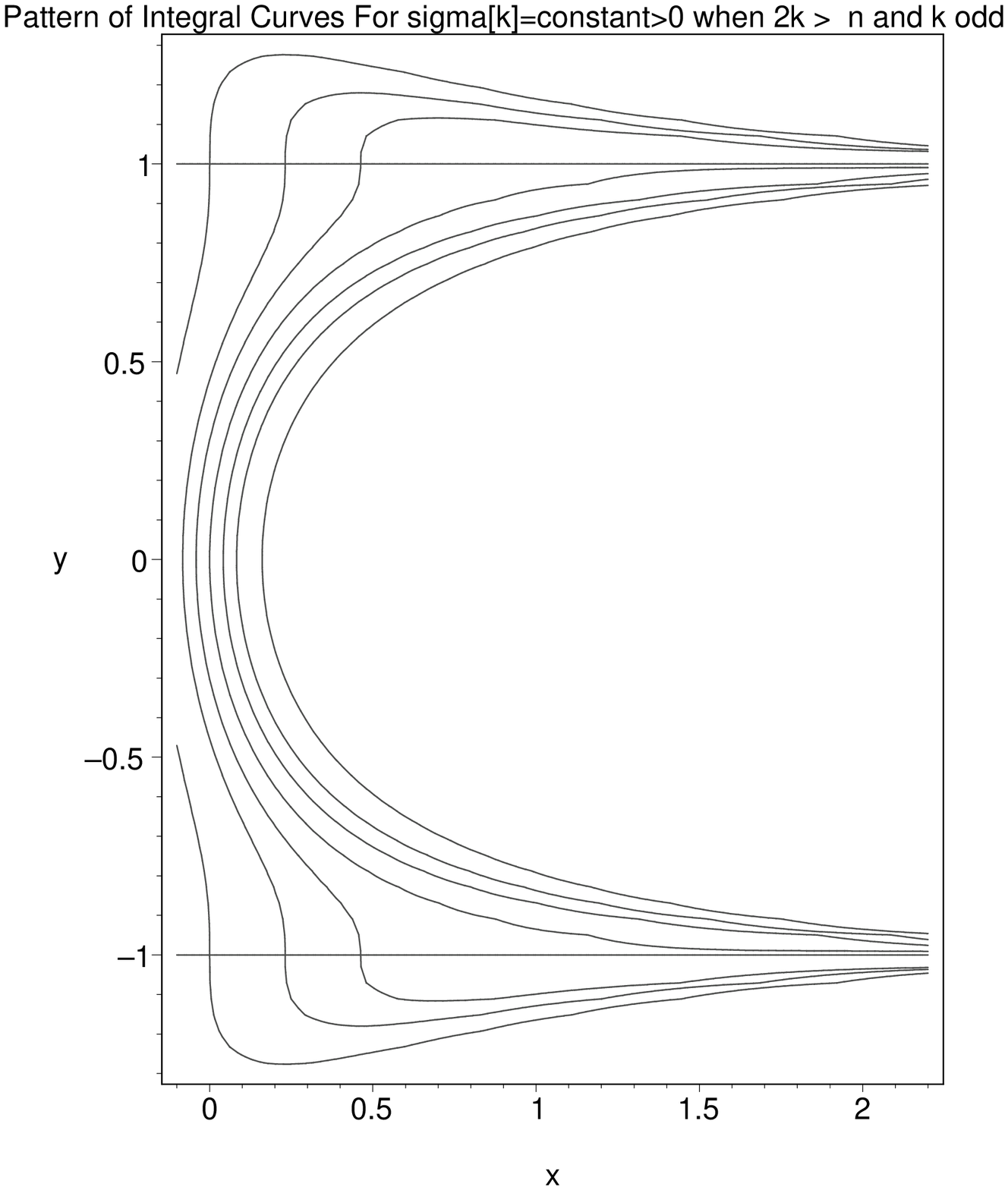}}\label{2kgnoddk}
\end{figure}
Despite the elementary nature of the analysis involved, the catalog of different behaviors of the radial
solutions is of potential reference value, so we provide a full catalog of the
 behaviors of the radial solutions in the following three theorems. 
For any $C^2$ solution of \eqref{3},
let $-\infty \leq t_{-} < t < t_{+}\leq \infty$ be
its maximum domain of definition, and correspondingly 
 $r_{-}=e^{t_{-}}$, and $r_{+}=e^{t_{+}}$.  Note that \eqref{3}, 
when $\sigma_k$ is a constant, is invariant
under the reflection $t \mapsto -t$, which, in terms of $|x|$, corresponds to
symmetry under inversion. The  results in the following  theorems are stated subject to this inversion.
\begin{theorem} 
The behavior of the radial solutions of \eqref{3}, when $k>1$ and 
 $\sigma_k$ is a \emph{positive} constant, normalized to be $2^{-k}\binom{n}{k}$, 
are cataloged as follows. Recall that, along any solution, either $1-\xi^2_t<0$, or $1-\xi^2_t>0$; and
$e^{(2k-n)\xi}(1-\xi_t^2)^k - e^{-n\xi}$ is a contant. Denote this constant by $h$.
\begin{description}
\item[Case I.] $1-\xi_t^2 > 0$. Recall that all such solutions are in the $\Gamma^{+}_{k}$ class.
These solutions fall into one of the
    following three categories.
  \begin{enumerate}
	\item  If $h=0$, then the domain of definition of $v(|x|)$ is the entire $\mathbb
	  R^n$, and $v(|x|)^{-2} = \left(\frac{2\rho}{|x|^2+\rho^2}\right)^2$ for some positive parameter
	$\rho$. So these solutions  give rise to the round
	  spherical metric on $\mathbb R^n \cup \{ \infty \}= \mathbb S^n$.
	\item If $h<0$, then the domain of definition of $v(|x|)$ is given by 
          $0 < r_{-} < |x| < r_{+} < \infty$. $v$ and $v_r$ stay bounded on
           $[r_{-}, r_{+}]$, in fact, $v_r \to 0$ as $r \to r_{-}$ and 
            $rv_r/v \to 2$ as $r \to r_{+}$, but $v_{rr}$ blows up at both ends of the
	  interval $[r_{-}, r_{+}]$.
	\item If $h>0$, then the behavior of $v$ is classified according to the relation between 
	      $2k$ and $n$:
		\begin{enumerate} 
			\item If $2k < n$, then $h$ has the further restriction
			$h \leq  h^*=\frac{2k}{n-2k}\left( \frac{n-2k}{n} \right)^\frac{n}{2k}$ and
			the domain of definition of $v(|x|)$ is given by 
          			$0 < |x| < \infty$. In fact,  $\xi (t)$ is a periodic
          			function of $t$, giving rise to a metric 
				$g = \frac{e^{-2\xi(\ln |x|)}}{|x|^2} |dx|^2 $ on
          			$\mathbb R^n \setminus \{0\}$ which is complete. Note that the case $h=h^*$
				gives rise to the cylindrical metric $\frac{|dx|^2}{|x|^2}$ on 
				$\mathbb R^n \setminus \{ 0\}$.
			\item If  $2k = n$, then $h$ satisfies the further restriction $h<1$ and
                                 the domain of definition of $v(|x|)$ is given by 
          			$0 < |x| < \infty$. As $|x| \to 0$, $v^{-2}(|x|)|dx|^2$ has the asymptotic
				$g \sim |x|^{-2(1- \sqrt{1-\sqrt[k]{h}})} |dx|^2$, and as $|x| \to \infty$,
				$v^{-2}(|x|)|dx|^2$ has the asymptotic 
				$g \sim |x|^{-2(1+ \sqrt{1-\sqrt[k]{h}})} |dx|^2$. Thus  
				$g$ gives rise to a metric on  $\mathbb R^n
          			\setminus \{0\}$ singular at $0$ and at $\infty$ which
          			behaves like the cone metric, is incomplete with
          			finite volume.
			\item If $2k > n$, then the domain of definition of $v(|x|)$ is given by 
          			$0 < |x| < \infty$.  $v^{-2}(|x|)$ has an asymptotic 
				expansion of the form 
		\[
		v^{-2}(|x|) = \rho^{-2}\{1- \sqrt[k]{h} \frac{k}{2k-n} 
				\left( \frac{|x|}{\rho}\right)^{2-\frac{n}{k}}+ \cdots \}
		\] 
				as $|x| \to 0$, where $\rho>0$ is a positive parameter,
				thus $v(|x|)$ has a positive, finite limit,  but $v_{rr}(|x|)$
				 blows up at $|x| \to 0$. The behavior of $v$
				as $|x|\to \infty$ can be described similarly. Putting together, 
				we conclude that
				$v^{-2}|dx|^2$ extends to a $C^{2-\frac{n}{k}}$ metric on $\mathbb S^n$.
                \end{enumerate}
  \end{enumerate}
\item[Case II.] $1-\xi_t^2 <0$ and $k$ even. Recall that all such solutions are in the 
$\Gamma^{-}_{k}$ class. Subject to an inversion, these solutions are defined for
$0\leq r_{-} < (\,  \text{or} \le \,) \, |x| < r_{+} < \infty$, and has the asymptotic expansion
$v^{-2} \sim (r_{+}-r)^{-2}$ as $|x| \to r_{+}$. Their behavior as $|x| \to r_{-}$ falls into one of the
    		following three categories.
  \begin{enumerate}
	 \item If $h=0$, then the  domain of definition of $v(|x|)$ is $|x| < r_+<
	       \infty$. These solutions define  the hyperbolic metric
	       defined on $|x| < r_+$.
	  \item If $h<0$, then the  domain of definition of $v(|x|)$ is 
                 $0< r_- < |x| < r_+< \infty$. As $|x| \to r_{-}$, $v(|x|)$ has a positive, finite
		limit, $v_r(|x|) \to 0$, but $v_{rr}(|x|)$ blows up. These solutions define
	         metrics on  $ r_- < |x| < r_+$, which is complete 
                  near $|x| = r_+$, and has its second derivative
                 blowing up as $|x| \to r_-$.
	  \item If $h>0$, then the behavior of $v$ as $|x|\to r_-$ is  classified according to the relation between 
	        $2k$ and $n$:
		\begin{enumerate}
			\item If $2k<n$, then the  domain of definition of $v(|x|)$ is 
	         		$0< r_- < |x| < r_+< \infty$. The corresponding metric  has
	        		 the following degeneracy at $r_-$: 
		 		$g \sim (r-r_-)^{ \frac {4k}{n-2k}} |dx|^2$ as $|x| \to r_-$, 
				and is complete as $|x| \to r_+$.
	  		\item If $2k=n$, then the  domain of definition of $v(|x|)$ is
	     			$0 < |x| < r_+ < \infty $. The metric has  the conical degeneracy 
	     			$g \sim |x|^{2(\sqrt{1+\sqrt[k]{h}}-1)} |dx|^2$ as $|x| \to 0$,
	     			and is complete as $|x| \to r_+$.
	 		\item If $2k>n$, then the  domain of definition of $v(|x|)$ is
	         		$0< |x| < r_+$.
				$v^{-2}(|x|)$ has an asymptotic 
				expansion of the form 
		\[
		v^{-2}(|x|) = \rho^{-2}\{1+ \sqrt[k]{h} \frac{k}{2k-n} 
				\left( \frac{|x|}{\rho}\right)^{2-\frac{n}{k}}+ \cdots \}
		\] 
				as $|x| \to 0$, where $\rho>0$ is a positive parameter,
				thus,  as $|x| \to 0$,  $v(|x|)$ has a positive, finite limit,  but $v_{rr}(|x|)$
				 blows up. $v^{-2}(|x|)|dx|^2$ is complete as $|x| \to r_+$.
	        \end{enumerate}
  \end{enumerate}
\item[Case III.] $1-\xi_t^2 <0$ and $k$ odd. In this case $h<0$. Subject to an inversion, 
these solutions are defined for
$0\leq r_{-} < (\,  \text{or} \le \,)\, |x| < r_{+} < \infty$, and as $r \to r_{+}$, $v$ and $v_r$
			stay bounded, but $v_{rr}$ blows up. The behavior of $v$ as $|x|\to r_-$ is  classified
 according to the relation between  $2k$ and $n$:
	\begin{enumerate}
	 \item If $2k<n$, then the  domain of definition of $v(r)$ is  $0 < r_{-} < r < r_{+} < \infty$. As 
			$r \to r_{-}$, $v^{-2}|dx|^2$ has the degeneracy: 
			$g \sim (r-r_{-})^{\frac{4k}{n-2k}}dx|^2$. 
	 \item If $2k = n$, then the  domain of definition of $v(|x|)$ is
	     			$0 < |x| < r_+ < \infty$. As $|x| \to 0$, $v^{-2}|dx|^2$ has the conical degeneracy 
	     			$g \sim |x|^{2(\sqrt{1+\sqrt[k]{|h|}}-1)} |dx|^2$.
	 \item If $2k> n$, then the  domain of definition of $v(|x|)$ is
	     			$0 < |x| < r_+ < \infty$. As $|x| \to 0$, $v^{-2}$ has the asymptotic
			expansion of the form 
		\[
		v^{-2}(|x|) = \rho^{-2}\{1+ \sqrt[k]{|h|} \frac{k}{2k-n} 
				\left( \frac{|x|}{\rho}\right)^{2-\frac{n}{k}}+ \cdots \}
		\] 
				as $|x| \to 0$, where $\rho>0$ is a positive parameter,
				thus $v(|x|)$ stays bounded, but $v_{rr}(|x|)$
				 blows up both as  $|x| \to 0$ and as $|x| \to r_{+}$.

	\end{enumerate}
\end{description}
\end{theorem}

\begin{remark}
At an end point $0< r_* < \infty$ of the domain of definition of $v(r)$ which corresponds to $v$ 
having a positive, finite limit, the proof for Theorem 1 will give the rate of blowing up of
$v_{rr}$ as proportional to $|r-r_*|^{-1+\frac{1}{k}}$.
\end{remark}

The behavior of the radial solutions of \eqref{3}, when $k>1$ and 
 $\sigma_k$ is a \emph{negative} constant, normalized to be $-2^{-k}\binom{n}{k}$, 
is cataloged in the following theorem, and the phase plane pattern in the $\xi$-$\xi_t$ plane
is displayed on the next page.
\begin{figure}[hp!]
\scalebox{0.40}{\includegraphics{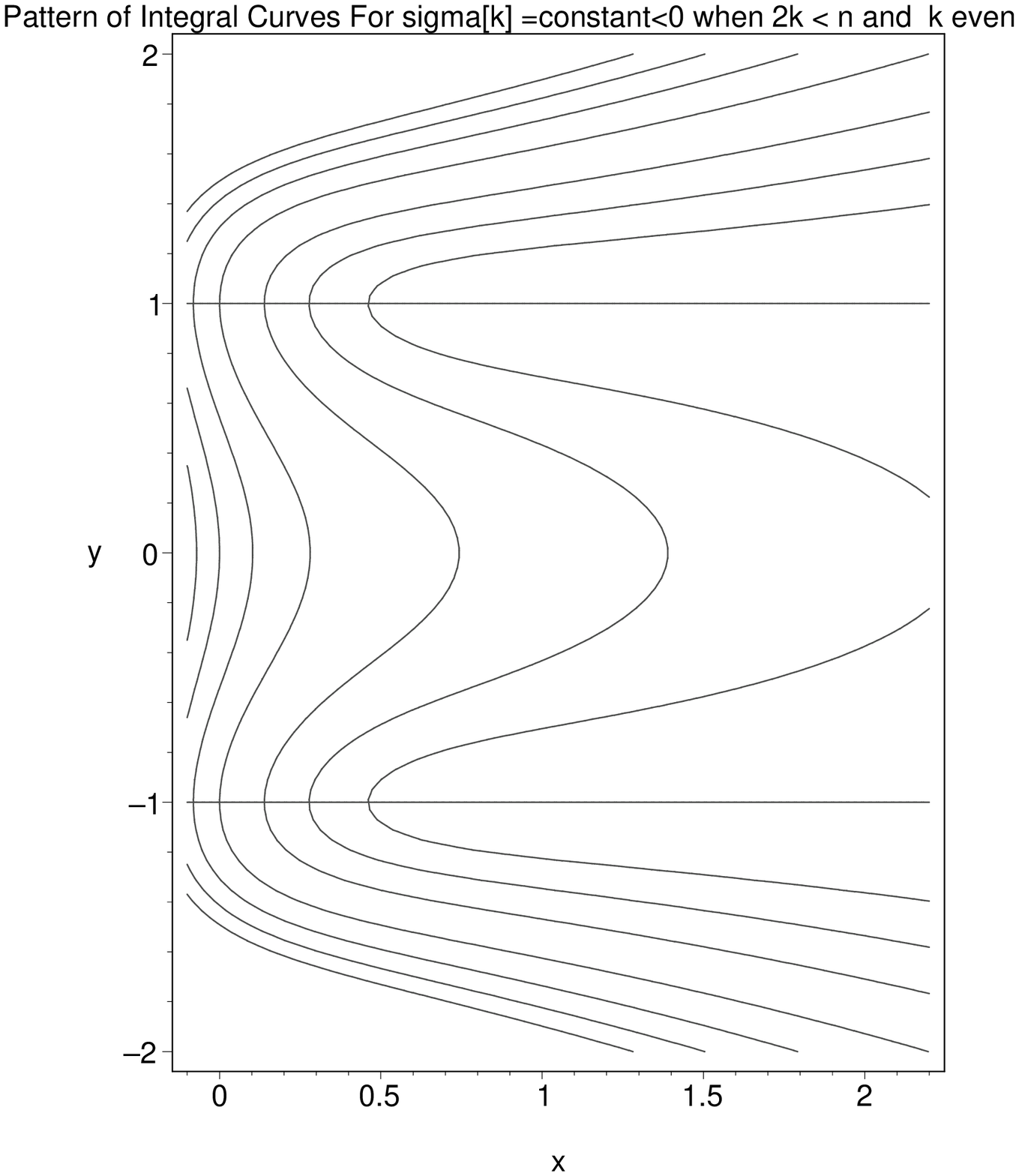}}\label{n2ksnevenk}
\scalebox{0.40}{\includegraphics{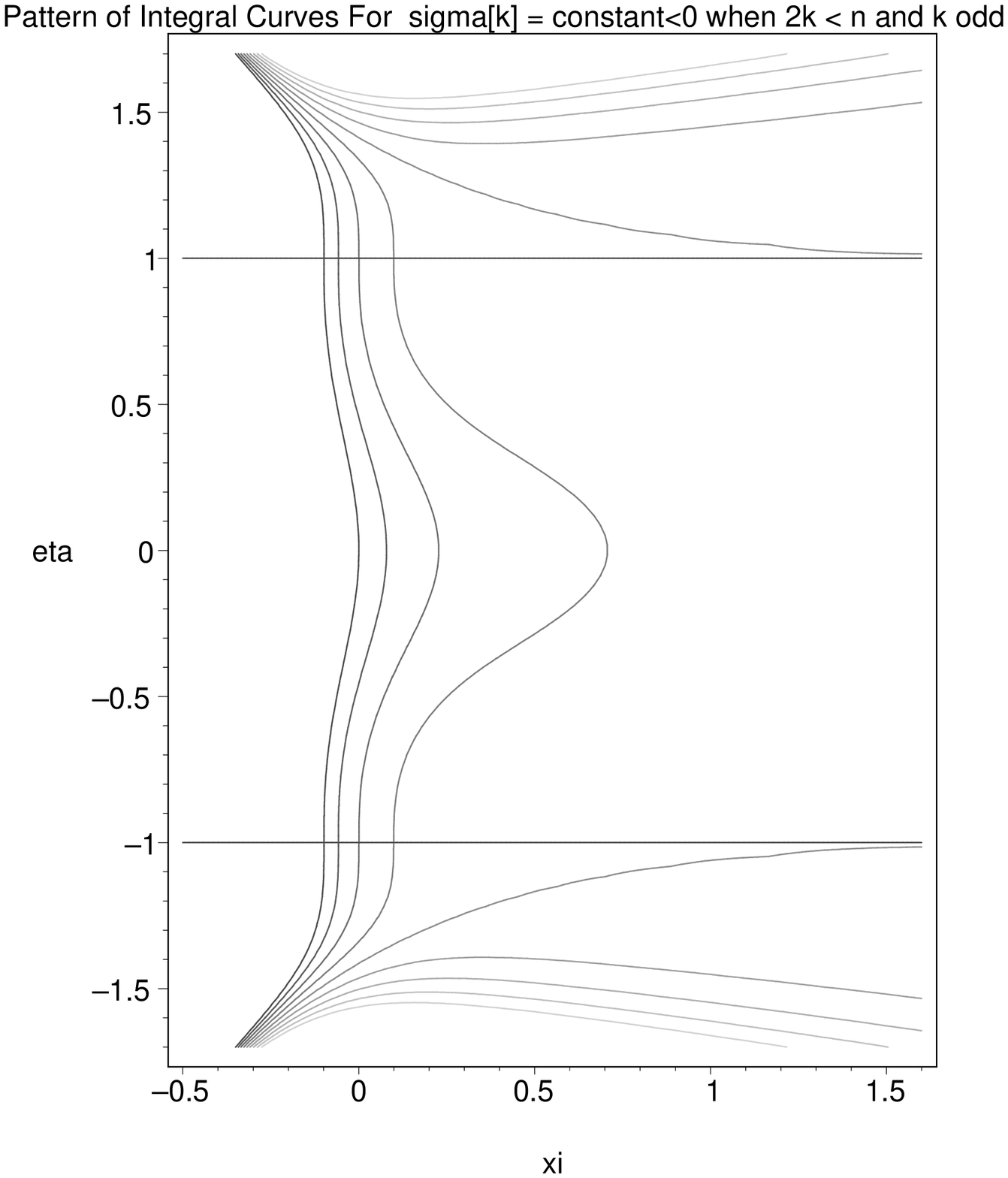}}\label{n2ksnoddk}
\scalebox{0.40}{\includegraphics{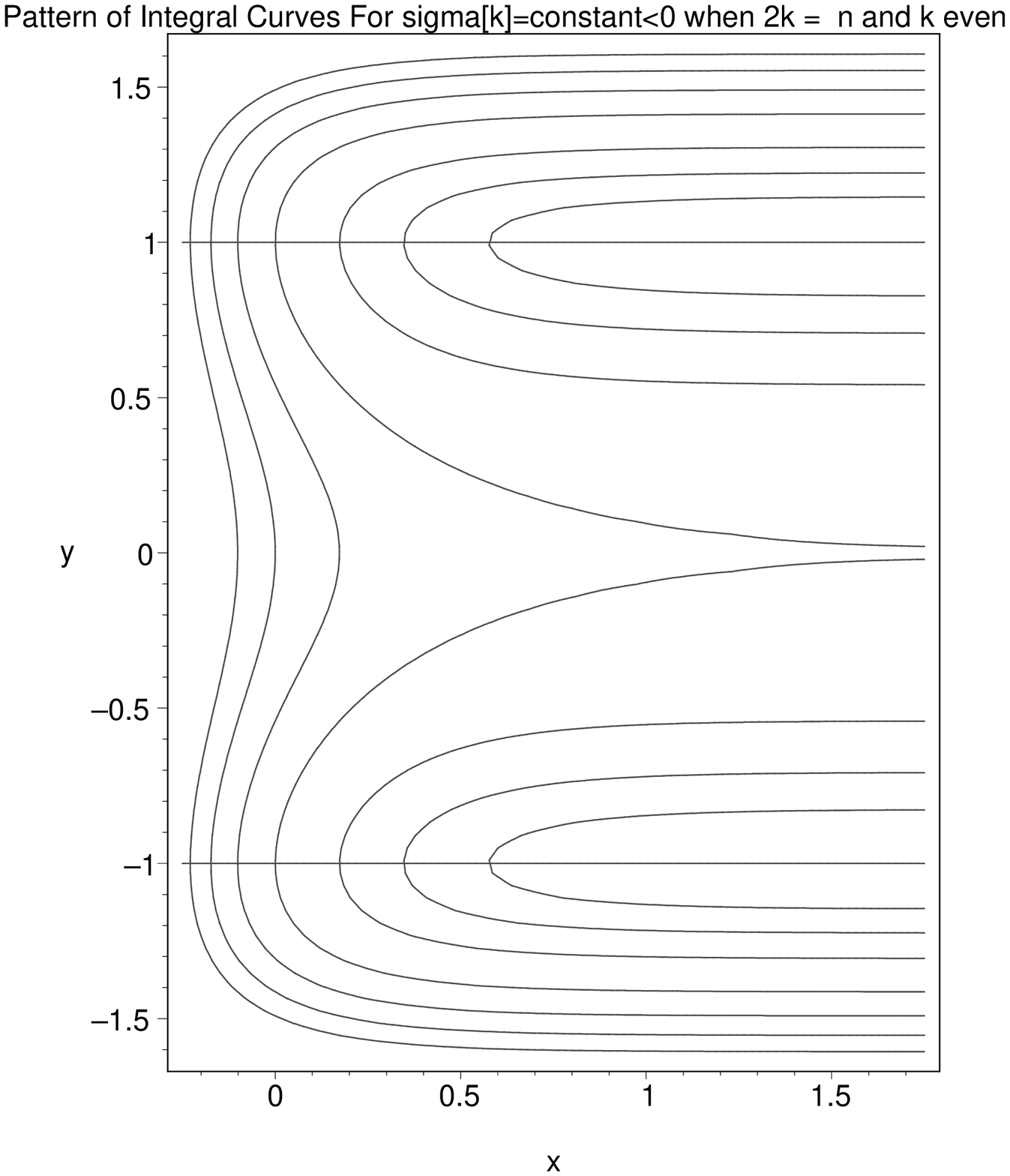}}\label{n2kenevenk}
\scalebox{0.40}{\includegraphics{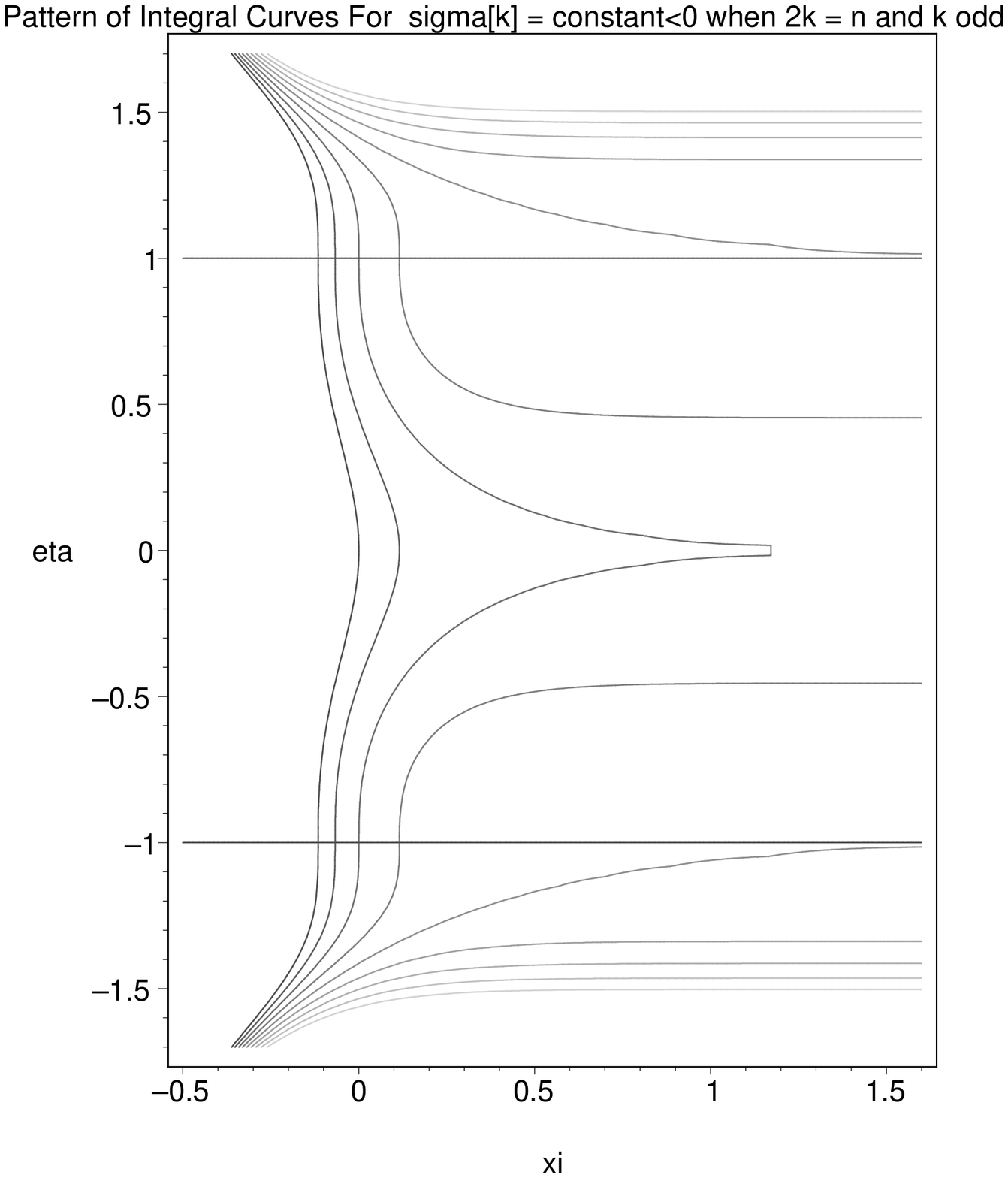}}\label{n2kenoddk}
\scalebox{0.40}{\includegraphics{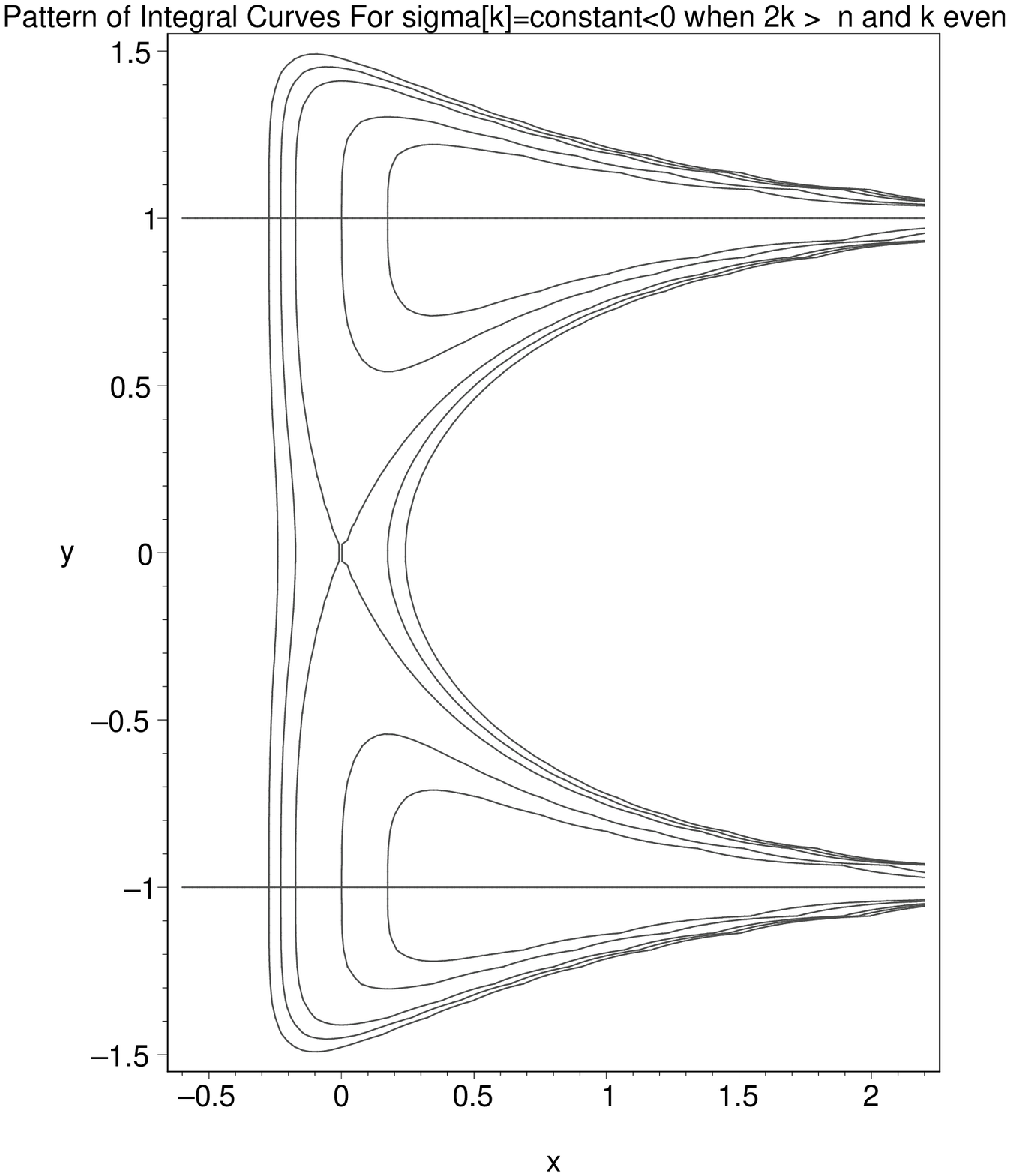}}\label{n2kgnevenk}
\quad \qquad \qquad
\scalebox{0.40}{\includegraphics{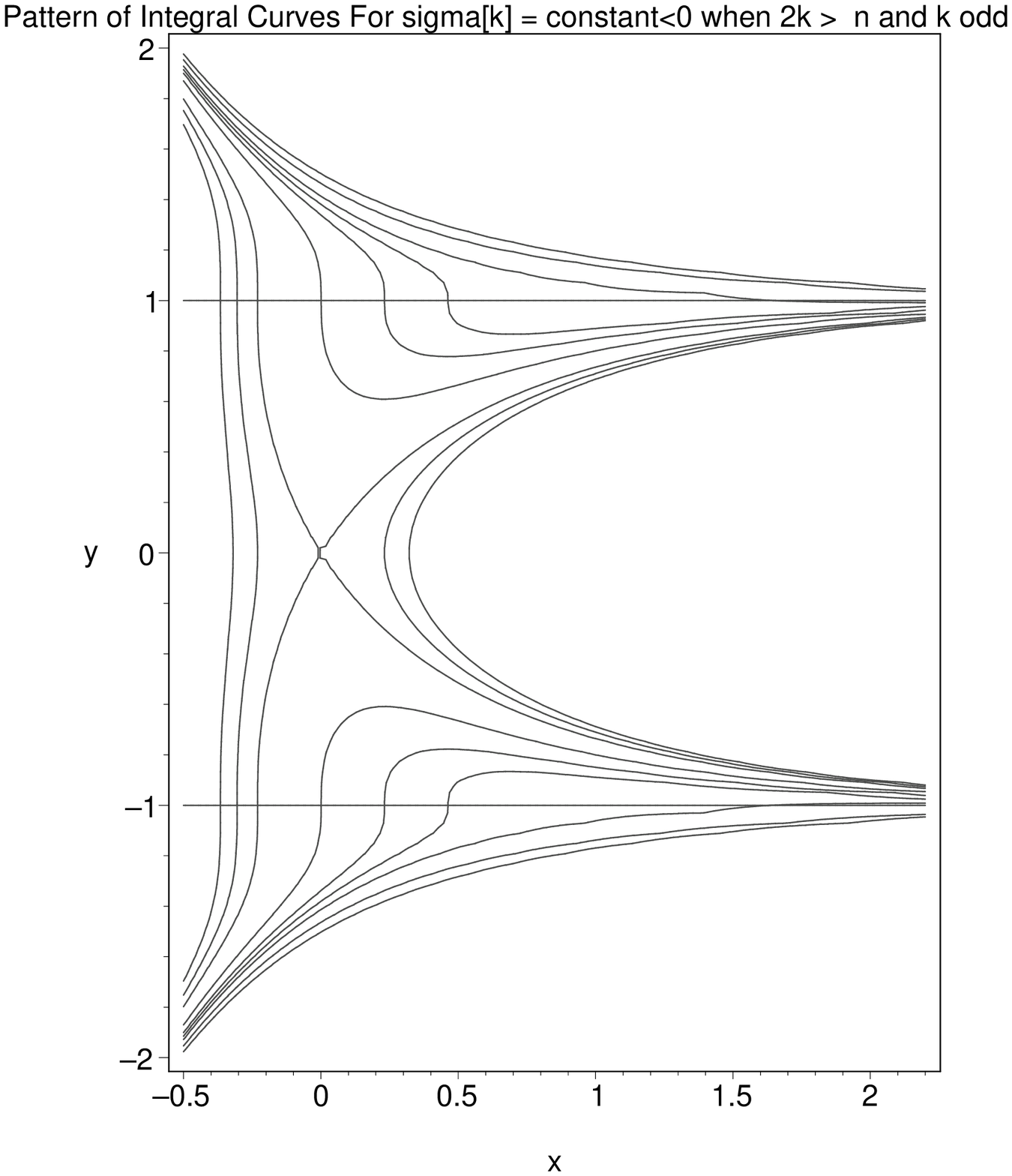}}\label{n2kgnoddk}
\end{figure}

\begin{theorem}  Recall that, along any solution of \eqref{3}, when
$k>1$ and $\sigma_k$ is a \emph{negative} constant, normalized to be $-2^{-k}\binom{n}{k}$, 
we have that either $1-\xi^2_t<0$, or $1-\xi^2_t>0$; and
$e^{(2k-n)\xi}(1-\xi_t^2)^k + e^{-n\xi}$ is a contant. Denote this constant by $h$.
\begin{description}
\item[Case I.]  $1-\xi_t^2 > 0$. All such solutions have $h>0$ and fall into one of the
    following three categories.
  \begin{enumerate}
	\item  If $2k <n$, then the domain of definition of $v(|x|)$ is given by 
          $0 < r_{-} < |x| < r_{+} < \infty$. The second
	  derivatives of these solutions  blow up at both ends of the
	  interval $[r_{-}, r_{+}]$.
	\item If $2k=n$, then $v$ or its inversion falls into one of the following three subcases:
		\begin{enumerate}
			\item If $0<h < 1$, then the domain of definition of $v(|x|)$ is given by
				$0 <|x| <r_{+}$. As $|x| \to r_{+}$, $v(|x|)$ has a positive, finite limit,
				$v_r(|x|) \to 0$, but $v_{rr}(|x|)$ blows up; as $|x| \to 0$, $v^{-2}(|x|)$
				has the asymptotic $v^{-2}(|x|) \sim |x|^{-2(1-\sqrt{1-\sqrt[k]{h}})}$. So the 
				metric $v^{-2}|dx|^2$ behaves like an incomplete, conic metric near $0$.
			\item If $h=1$, then the domain of definition of $v(|x|)$ is given by
				$0 <|x| <r_{+}$. As $|x| \to r_{+}$, $v(|x|)$ has a positive, finite limit,
				$v_r(|x|) \to 0$, but $v_{rr}(|x|)$ blows up; as $|x| \to 0$, $v^{-2}(|x|)$
				has the asymptotic $v^{-2}(|x|) \sim |x|^{-2}(\ln\frac{1}{|x|})^{-\frac{2}{k}}$.
				So the 
				metric $v^{-2}|dx|^2$ can be thought of as a complete metric near $0$.
			\item If $h>1$, then the domain of definition of $v(|x|)$ is given by 
          $0 < r_{-} < |x| < r_{+} < \infty$. The second
	  derivatives of these solutions  blow up at both ends of the
	  interval $[r_{-}, r_{+}]$.
		\end{enumerate}
	\item If $2k >n$, then $v$ or its inversion falls into one of the following  six  subcases. Let 
		$h^*=\frac{2k}{2k-n}\left( \frac{2k-n}{n} \right)^\frac{n}{2k}$. Note that 
	$M(h)=\max_{\xi \in \mathbb R} \{e^{(n-2k)\xi}h- e^{-2k\xi}\}$ is strictly monotone in
		$h$, and $h^*$ is the unique positive number such that $M(h^*)=1$.
		\begin{enumerate} 
			\item If $h < h^*$, then the domain of definition of $v(|x|)$ is given by 
          			$0 < |x| < r_{+}$. As $|x| \to r_{+}$, $v(|x|)$ has a positive, finite limit,
				$v_r(|x|) \to 0$, but $v_{rr}(|x|)$ blows up; as $|x| \to 0$, $v^{-2}(|x|)$
				has the asymptotic 
			\[
		v^{-2}(|x|) = \rho^{-2}\{1-  \sqrt[k]{h} \frac{k}{2k-n} 
				\left(\frac{|x|}{\rho}\right)^{2-\frac{n}{k}}+ \cdots \},
		\] 
				where $\rho>0$ is a positive parameter. This case corresponds to 
			     the integral curves in the phase portrait that intercept and are also asymptotic
				to the lines $\xi_t = \pm 1$.
			\item If $h = h^*$ and $\xi_{tt} <0$, then the domain of definition of $v(|x|)$ is given by 
          			$0 < |x| < r_{+}$. As $|x| \to r_{+}$, $v(|x|)$ has a positive, finite limit,
				$v_r(|x|) \to 0$, but $v_{rr}(|x|)$ blows up; as $|x| \to 0$, $v^{-2}(|x|)|dx|^2$
				has the cylindrical metric $|x|^{-2} |dx|^2$ as the asymptotic. This case
				corresponds to the integral curves in the phase portrait that are asymptotic
				to the unique equilibrium point 
				$\xi = \xi^* = \frac{1}{n} \ln [\frac{2k}{(2k-n)h}], \xi_t =0$ and intercept the
				lines $\xi_t = \pm 1$.
 			\item If $h = h^*$ and $\xi_{tt} >0$, then the domain of definition of $v(|x|)$ is given by
				$0 < |x| < \infty$. As $|x| \to \infty$, $v^{-2}(|x|)|dx|^2$
				has the cylindrical metric $|x|^{-2} |dx|^2$ as the asymptotic. As $|x| \to 0$,
				$v^{-2}(|x|)$
				has the asymptotic 
			\[
		v^{-2}(|x|) = \rho^{-2}\{1-  \sqrt[k]{h} \frac{k}{2k-n} 
				\left(\frac{|x|}{\rho}\right)^{2-\frac{n}{k}}+ \cdots \},
		\]
			where $\rho>0$ is a positive parameter. This case
				corresponds to the integral curves in the phase portrait that are asymptotic
				to both the unique equilibrium point 
				$\xi = \xi^* = \frac{1}{n} \ln [\frac{2k}{(2k-n)h}], \xi_t =0$ and  the
				lines $\xi_t = \pm 1$.
			\item If $h = h^*$ and $\xi_{tt} \equiv 0$, then $\xi \equiv \xi^*$ and the domain of 
				definition of $v(|x|)$ is given by
				$0 < |x| < \infty$. The metric  $v^{-2}(|x|)|dx|^2$ is the cylindrical metric 
				$|x|^{-2} |dx|^2$.
			\item If  $h>h^*$ and $\xi_{tt}<0$, then the domain of definition of $v(|x|)$ is given by 
          			$0 < r_{-} < |x| < r_{+} < \infty$. The second
	  			derivatives of these solutions  blow up at both ends of the
	  			interval $[r_{-}, r_{+}]$. This case corresponds to the integral curves in the 
				phase portrait that intercepts the lines $\xi_t = \pm 1$ at both ends.
			\item If $h>h^*$ but $\xi_{tt}>0$, then the domain of definition of $v(|x|)$ is given by 
          			$0 < |x| < \infty$, and as $|x| \to 0$, $v^{-2}(|x|)$ has the same asymptotic
				as in case (a) above; as $|x| \to \infty$, $v^{-2}(|x|)$ has a similar asymptotic.
				One may think of $v^{-2}(|x|)|dx|^2$ as extending to a $C^{2-\frac{n}{k}}$ metric
				on $\mathbb S^n$. This case corresponds to the integral curves in the 
				phase portrait that are asymptotic to the lines $\xi_t = \pm 1$ at both ends.
                \end{enumerate}
  \end{enumerate}
\item[Case II.] $k$ odd and $1-\xi_t^2 <0$. Recall that all such solutions are in the 
$\Gamma^{-}_{k}$ class. Subject to an inversion, these solutions are defined for
$0\leq r_{-} < (\, \text{or} \le \,)\, |x| < r_{+} < \infty$, and has the asymptotic expansion
$v^{-2} \sim (r_{+}-r)^{-2}$ as $|x| \to r_{+}$. Their behavior as $|x| \to r_{-}$ falls into one of the
    		following three categories.
  \begin{enumerate}
	 \item If $h=0$, then the  domain of definition of $v(|x|)$ is $|x| < r_+<
	       \infty$. These solutions define  the hyperbolic metrics
	       defined on $|x| < r_+$.
	  \item If $h>0$, then the  domain of definition of $v(|x|)$ is 
                 $0< r_- < |x| < r_+< \infty$. As $|x| \to r_-$, $v(|x|)$ has a positive limit, 
			$v_r(|x|) \to 0$, but $v_{rr}(|x|)$ blows up. 
	  \item If $h<0$, then the behavior of $v$ as $|x| \to r_{-}$
		is classified according to the relation between $2k$ and $n$:
		\begin{enumerate}
			\item If $2k<n$, then the  domain of definition of $v(|x|)$ is 
	         		$0< r_- < |x| < r_+< \infty$. The metric $v^{-2}|dx|^2$ has
	        		 the degeneracy at $r_-$: 
		 		$g \sim (r-r_-)^{ \frac {4k}{n-2k}} |dx|^2$ as $|x| \to r_-$, 
				and is complete as $|x| \to r_+$.
	  		\item If $2k=n$, then the  domain of definition of $v(|x|)$ is
	     			$0 < |x| < r_+ < \infty$. The metric $v^{-2}|dx|^2$ has the conical degeneracy 
	     			 $g \sim |x|^{2(\sqrt{1+\sqrt[k]{|h|}}-1)} |dx|^2$ as $|x| \to 0$,
	     			and is complete as $|x| \to r_+$.
	 		\item If $2k>n$, then  the  domain of definition of $v(|x|)$ is
	     			$0 < |x| < r_+ < \infty$. As $|x| \to 0$, $v^{-2}$ has the asymptotic
			expansion of the form 
		\[
		v^{-2}(|x|) = \rho^{-2}\{1+ \sqrt[k]{|h|} \frac{k}{2k-n} 
				\left( \frac{|x|}{\rho}\right)^{2-\frac{n}{k}}+ \cdots \}
		\] 
				as $|x| \to 0$, where $\rho>0$ is a positive parameter.
	        \end{enumerate}
  \end{enumerate}
\item[Case III.] $k$ even and $1-\xi_t^2 <0$. In this case $h>0$. All such solutions, or their inversion,
	 are defined for $0\le r_{-} < (\text{or} \le \,)\, |x| < r_{+} < \infty$, and as $r \to r_{+}$, $v$ and 
	$v_r$ stay bounded, but $v_{rr}$ blows up. The behavior of $v$ as $|x| \to r_{-}$
		is classified according to the relation between $2k$ and $n$:
		
	\begin{enumerate}
	 \item If $2k<n$, then the  domain of definition of $v(r)$ is  $0 < r_{-} < r < r_{+} < \infty$. As 
			$r \to r_{-}$, $v^{-2}|dx|^2$ has the degeneracy: 
			$g \sim (r-r_{-})^{\frac{4k}{n-2k}} |dx|^2$.
	 \item If $2k = n$, then the  domain of definition of $v(|x|)$ is
	     			$0 < |x| < r_+ < \infty$. As $|x| \to 0$, $v^{-2}|dx|^2$ has the conical degeneracy 
	     			$g \sim |x|^{2(\sqrt{1+\sqrt[k]{h}}-1)} |dx|^2$.
	 \item If $2k> n$, then the  domain of definition of $v(|x|)$ is
	     			$0 < |x| < r_+ < \infty$. As $|x| \to 0$, $v^{-2}$ has the asymptotic
			expansion of the form 
		\[
		v^{-2}(|x|) = \rho^{-2}\{1+ \sqrt[k]{h} \frac{k}{2k-n} 
				\left( \frac{|x|}{\rho}\right)^{2-\frac{n}{k}}+ \cdots \}
		\] 
				as $|x| \to 0$, where $\rho>0$ is a positive parameter,
				thus $v(|x|)$ stays bounded, but $v_{rr}(|x|)$
				 blows up both as  $|x| \to 0$ and as $|x| \to r_{+}$.

	\end{enumerate}

\end{description}
\end{theorem}
The classification of solutions of \eqref{3} when $\sigma_k(A_g) \equiv 0$ is easily done by
integrating out the equation.
\begin{theorem}
Any radial solution of \eqref{3}, when $\sigma_k(A_g) \equiv 0$, is one of the following three forms.
\begin{enumerate}
\item $\xi \equiv \pm 1$;
\item $\xi (t) = (1 -\frac{n}{2k})^{-1} \ln \left|\sinh \left[(1 -\frac{n}{2k})(t-t_0)\right] \right| + c$ for some
constants $t_0, c$;
\item $\xi (t) = (1 -\frac{n}{2k})^{-1} \ln \cosh\left[(1 -\frac{n}{2k})(t-t_0)\right] +c$
for some constants $t_0, c$.
\end{enumerate}
\end{theorem}

\section{Indication of Proofs.}
Here is an elementary proof of Proposition 1.
\begin{proof}  First, the conclusion in part 1 follows readily from \eqref{3}.
To prove part 2, notice that $\sigma_k >0$ and $1-\xi_t^2>0$ imply that 
$\frac{k}{n} \xi_{tt} + (\frac{1}{2} -\frac{k}{n})(1-\xi_t^2) >0$. Then it follows, for $1\leq l<k$, that
\begin{eqnarray*}
\sigma_l(A_g) &=& c_{n,l}^{'}(1-\xi_t^2)^{l-1}\left[\frac{l}{n} \xi_{tt}+(\frac{1}{2}-\frac{l}{n})(1-\xi_t^2)\right] e^{2l\xi} \\
        &=& \frac{l}{k} c_{n,l}^{'} (1-\xi_t^2)^{l-1} \left[\frac{k}{n} \xi_{tt} + (\frac{k}{2l} -\frac{k}{n})(1-\xi_t^2)\right] e^{2l\xi}\\
              &= &\frac{l}{k} c_{n,l}^{'} (1-\xi_t^2)^{l-1} \left[\frac{k}{n} \xi_{tt} + (\frac{1}{2} -\frac{k}{n})(1-\xi_t^2)+
                  (\frac{k}{2l} -\frac{1}{2})(1-\xi_t^2)\right] e^{2l\xi} \\
	      &>& 0.
\end{eqnarray*}
Similarly, in the situation for part 3, $\sigma_k >0$, $1-\xi_t^2<0$, and $k$ is even, then
$\frac{k}{n} \xi_{tt} + (\frac{1}{2} -\frac{k}{n})(1-\xi_t^2)  <0$, and
\begin{eqnarray*}
(-1)^l\sigma_l(A_g) &=&(-1)^l\frac{l}{k} c_{n,l}^{'} (1-\xi_t^2)^{l-1} \left[\frac{k}{n} \xi_{tt} + 
(\frac{1}{2} -\frac{k}{n})(1-\xi_t^2)+ (\frac{k}{2l} -\frac{1}{2})(1-\xi_t^2)\right] e^{2l\xi} \\
              &>& 0.
\end{eqnarray*}
Finally, in the situation for part 4, $\sigma_k < 0$, $1-\xi_t^2<0$, and $k$ is odd, then
$\frac{k}{n} \xi_{tt} + (\frac{1}{2} -\frac{k}{n})(1-\xi_t^2) <0$, and
\begin{eqnarray*}
(-1)^l\sigma_l(A_g) &=&(-1)^l\frac{l}{k} c_{n,l}^{'} (1-\xi_t^2)^{l-1} \left[\frac{k}{n} \xi_{tt} + 
(\frac{1}{2} -\frac{k}{n})(1-\xi_t^2)+(\frac{k}{2l} -\frac{1}{2})(1-\xi_t^2)\right] e^{2l\xi} \\
              &>& 0.
\end{eqnarray*}
For part 5, recall that the equation in this case has a first integral: when $\sigma_k(A_g)$ is
normalized to be $\pm 2^{-k} \binom{n}{k}$, in terms of $\xi = \ln (v/r)$ we have that
$e^{(2k-n)\xi}(1-\xi_t^2)^k \pm e^{-n\xi}$
is a constant along any solution. The conclusion of part 5 follows then easily, noting the 
relation $1-\xi_t^2 = (2-rv_r/v)rv_r/v$.
\end{proof}
For the proof to Theorems 1 and 2, we will offer enough details for the case $\sigma_k = 2^{-k} \binom{n}{k}$,
 but leave out
the details for the case of $\sigma_k = - 2^{-k} \binom{n}{k}$.

\begin{proof}[Proof of Theorem 1]  Recall that, in this case,
 for each solution $\xi (t)$, there is a constant $h$ such that
\begin{equation} \label{h}
(1-\xi^2_t)^k = e^{-2k\xi} + h e^{(n-2k)\xi},
\end{equation}
along the solution. Set $D(\xi) = e^{-2k\xi} + h e^{(n-2k)\xi}$. Recall also that either $1-\xi_t^2 > 0$ or 
$1-\xi_t^2 < 0$ along the solution.

\noindent {\bf Case I.} \underline{$1-\xi_t^2 > 0$}. With the phase plane portrait as a guidance, 
it is routine to see that:
\begin{enumerate}
\item  If $h=0$, then $1 - \xi_t^2 = e^{-2\xi}$, which can be integrated out to produce 
$\xi = \ln \cosh (t-c)$. In terms of $v$, we see that $v^{-2} = (\frac{2e^c}{|x|^2+e^{2c}})^2$.
So $v^{-2}|dx|^2$ gives rise to the round spherical metric on $\mathbb S^n$ in this case.

\item If $h< 0$, then $\xi^{-} \le \xi < \xi^{+}$, where $\xi^{\pm}$ are (unique) finite numbers determined by
$D(\xi^{-}) =1$ and $D(\xi^{+}) = 0$. As $\xi \to \xi^{+}$, 
$\xi_t^2 \to 1$, which implies that $t$ tends to finite limits, thus the maximum domain of definition
of $v$ is $r_{-} <|x| < r_{+}$ for some $0< r_{-}  < r_{+}< \infty$, and because $\xi_t^2 \to 1$ 
and $\xi \to \xi^{+}$ as $|x| \to r_{\pm}$, $v$ and $v_r$ stay bounded, but $v_{rr}$ blows up as 
$|x| \to r_{\pm}$. Also as $r \to r_+$, $\xi_t \to 1$, thus $rv_r/v \to 2$; while as $r \to t_-$, $\xi_t \to -1$,
thus $v_r \to 0$.

\item If $h > 0$, then the analysis depends on the relation between $2k$ and $n$: 
\begin{enumerate}
	\item If $2k < n$, then $\xi^{-} \leq \xi \leq \xi^{+}$, where $\xi^{\pm}$ are the two roots of
$D(\xi)= 1$. This case requires that $h \leq h^*_{n,k}$, where $h^*_{n,k}$ is the unique positive
number such that $\min_{\mathbb R} D(\xi) =1$.  In this case,
the phase portrait shows that the solution curve $(\xi(t), \xi_t(t))$ is a periodic
orbit   bounded between $\xi^{-}$ and $\xi^{+}$, thus is defined for all $t$. So the metric 
$g = \frac{e^{-2\xi(\ln |x|)}}{|x|^2} |dx|^2$ is a complete metric on $\mathbb R^n \setminus \{0\}$.

	\item If $2k = n$, then $\xi^{-} \leq \xi$, where
$D(\xi^{-}) = 1$, and as $\xi \to \infty$, $D(\xi) \to h$, so $0< h < 1$ and 
$\xi_t^2 \to 1 - \sqrt[k]{h}$ as $\xi \to \infty$. Using $D(\xi) = h + e^{-2k\xi}$ and the relation
\begin{equation} \label{-k}
dt = \pm \frac {d\xi}{\sqrt{1-\sqrt[k]{D(\xi)}}},
\end{equation}
which follows from \eqref{h}, we conclude that 
\begin{enumerate}
\item[(i)] $\xi(t)$ is defined
for $-\infty < t < \infty$, and 
\item[(ii)] $t \pm \frac {\xi}{\sqrt{1-\sqrt[k]{h}}}$ has a finite limit
as $\xi \to \infty$.
\end{enumerate}
 Therefore, near $|x| \sim 0$, we have
\[
g \sim |x|^{-2(1- \sqrt{1-\sqrt[k]{h}})} |dx|^2,
\]
and near $|x| \sim \infty$, we have
\[
g \sim |x|^{-2(1+ \sqrt{1-\sqrt[k]{h}})} |dx|^2.
\]
These are incomplete, finite volume metrics on $\mathbb R^n \setminus \{0\}$, corresponding to
conical metrics on $\mathbb S^n \setminus \{0, \infty \}$.  

	\item If $2k > n$, then $\xi^{-} \leq \xi$, where
$D(\xi^{-}) = 1$, and as $\xi \to \infty$, $D(\xi) \to 0$ and $\xi_t^2 \to 1$.
Using \eqref{-k} and the asymptotic expansion $\sqrt[k]{D} \sim \sqrt[k]{h} e^{\frac{n-2k}{k}\xi}$,
we conclude that 
\begin{enumerate}
\item[(i)] $\xi(t)$ is defined
for $-\infty < t < \infty$, and 
\item[(ii)]  $ \xi \pm t = c + \frac{\sqrt[k]{h}}{2} \frac{k}{2k-n} e^{-\frac{2k-n}{k}\xi} + h.o.t.$
as $\xi \to \infty$,  
\end{enumerate}
for some constant $c$, from which we conclude that as $|x| \to 0$, we take the $+$ sign in (ii) and 
\[
v^{-2} = e^{-2c} \{1- \sqrt[k]{h} \frac{k}{2k-n} e^{-\frac{2k-n}{k}c}|x|^{\frac{2k-n}{k}}+ h.o.t. \}
\]
as $|x| \to 0$. The analysis near $x$ at $\infty$ can be carried out in a similar way. Thus we conclude
that $v^{-2}|dx|^2$ extends to  a $C^{2 - \frac{n}{k}}$ metric on $\mathbb S^n$ .
\end{enumerate}
\end{enumerate}

 The analysis for the
case $1-\xi_t^2 <0$ depends on whether $k$ is odd or even.

\noindent {\bf Case II.}\underline{ $1-\xi_t^2 <0$ and $k$ even}.  
\begin{enumerate}
\item If $h=0$, then $\xi_t^2 = 1 + e^{-2\xi}$. From this, we obtain
\[
dt = \pm \frac {d\xi}{\sqrt{1+e^{-2\xi}}}.
\]
It is easy to conclude from here that $v^{-2}|dx|^2$ provides the hyperbolic metric on either
$\{ |x| < r \}$ or $\{ |x| > r \}$.

\item If $h< 0$, then $\xi < \xi^{+}$, 
where $\xi^{+}$ is the
unique root of $D(\xi)=0$. As $\xi \to \xi^{+}$, $\xi_t^2 \to 1$ and $t$ has a finite limit. In terms of $v$ and
$|x|$, this produces a boundary point where $v_{rr}$ blows up. As $\xi \to -\infty$, $D(\xi) \to \infty$. 
In fact $\sqrt[k]{D(\xi)}  \sim e^{-2\xi}$.
Using this and  
\begin{equation} \label{+k}
dt = \pm \frac {d\xi}{\sqrt{1+\sqrt[k]{D(\xi)}}},
\end{equation}
we conclude that $t$ also has a finite limit. Taking the $-$ sign, for instance, and letting $t_{+}$ denote
this (upper) limit of $t$, we obtain
\[
t_{+} - t = e^{\xi} + h.o.t.
\]
as $\xi \to -\infty$, so that
\[
g \sim (t_{+}-t)^{-2} |dx|^2 = (\ln \frac{r_{+}}{|x|})^{-2} |dx|^2
\]
as $|x| \to r_{+}$. In conclusion, we obtain a conformal metric $v^{-2}|dx|^2$
on $r_{-} <|x| < r_{+}$ which has second derivative
blow up near $r_{-}$ and is complete near $r_{+}$.

\item If $h > 0$, the phase portrait indicates that the range for $\xi $ is the entire real line.
As $\xi \to -\infty$, $D(\xi) \to \infty$. In fact $\sqrt[k]{D(\xi)}  \sim e^{-2\xi}$.
Using this and \eqref{+k},
we conclude that $t$  has a finite limit as $\xi \to -\infty$, and again denoting the corresponding
(upper) limit of $|x|$ as $r_{+}$, we have
\[
g \sim (t_{+}-t)^{-2} |dx|^2 = (\ln \frac{r_{+}}{|x|})^{-2} |dx|^2,
\]
as $|x| \to r_{+}$. 
so $g$ is complete as $|x| \to r_{+}$. The analysis as $\xi \to \infty$ depends on the
relation between $2k$ and $n$.
	\begin{enumerate}
	\item If $2k < n$, then $\sqrt[k]{D(\xi)} \sim \sqrt[k]{h} e^{\frac{n-2k}{k} \xi}$ as $\xi \to \infty$.
Using this and  \eqref{+k},
we conclude that $t$  has a finite limit as $\xi \to \infty$. Taking the $-$ sign case, for instance,
and denoting this limit as $t_{-}$, we have
\[
t-t_{-} \sim e^{-\frac{n-2k}{2k}\xi}
\]
as $\xi \to \infty$. So
\[
g \sim (t-t_{-})^{\frac{4k}{n-2k}} |dx|^2 = (\ln \frac{|x|}{r_{-}})^{\frac{4k}{n-2k}} |dx|^2,
\]
and the metric $v^{-2}|dx|^2$ is defined on $r_{-}< |x| < r_{+}$, and as $|x| \to r_{-}$, it has the
above degeneracy; as $|x| \to r_{+}$, it is complete.

	\item If $2k = n$, then $\sqrt[k]{D(\xi)} \sim \sqrt[k]{h}$ as $\xi \to \infty$.
Using this and  \eqref{+k},
we conclude that $t \to \pm \infty$ as $\xi \to \infty$. In fact,
 $t \sim \pm \frac {\xi}{\sqrt{1+\sqrt[k]{h}}}$. Taking the $-$ sign, for instance, we obtain
\[
g \sim |x|^{2(\sqrt{1+\sqrt[k]{h}}-1)} |dx|^2.
\]
So the metric $v^{-2}|dx|^2$ is defined on $0 < |x| < r_{+}$, and as $|x| \to 0$, it has the
above degeneracy; as $|x| \to r_{+}$, it is complete.

	\item If $2k > n$, then $\sqrt[k]{D(\xi)} \to 0$ as $\xi \to \infty$. In fact,
$\sqrt[k]{D(\xi)} \sim \sqrt[k]{h} e^{\frac{n-2k}{k} \xi}$. 
Using this and  \eqref{+k},
we conclude that $t \to \pm \infty$ as $\xi \to \infty$. In fact, for some costant $c$,
\[
 \xi \pm t = c - \frac{\sqrt[k]{h}}{2} \frac{k}{2k-n} e^{-\frac{2k-n}{k}\xi} + h.o.t.
\]
as $\xi  \to \infty$, 
from which we conclude that
\[
v^{-2} = e^{-2c}\{1 + \sqrt[k]{h} \frac{k}{2k-n} e^{-\frac{2k-n}{k}c}|x|^{\frac{2k-n}{k}}+ h.o.t. \}
\]
as $|x| \to 0$. Thus the metric $v^{-2}|dx|^2$ defined on $0 < |x| < r_{+}$ is complete as $|x| \to r_{+}$,
and extends to be a $C^{2-\frac{n}{k}}$ metric on $|x| < r_{+}$.
	\end{enumerate}
\end{enumerate}

\noindent {\bf Case III.} \underline{ $1-\xi_t^2 <0$ and $k$ odd}. In this case $h<0$ and all such solutions have 
a finite limit point where $\xi_t = {\pm} 1$, which corresponds to a limit point $0<r_{*}<\infty$ where $v(|x|)$ has
a positive finite limite, but $v_{rr}$ blows up. Bu inversion, we may take $r_{+} = r_{*}$. The behavior of 
$v(|x|)$ as $|x| \to r_{-}$ depends on the
relation between $2k$ and $n$.
\begin{enumerate}
\item If $2k < n$, then $\sqrt[k]{D(\xi)} \sim \sqrt[k]{h} e^{\frac{n-2k}{k} \xi}$ as $\xi \to \infty$.
Using this and  \eqref{+k},
we conclude that $t$  has a finite limit as $\xi \to \infty$. Taking the $-$ sign case, for instance,
and denoting this limit as $t_{-}$, we have
\[
t-t_{-} \sim e^{-\frac{n-2k}{2k}\xi}
\]
as $\xi \to \infty$, so
\[
g \sim (t-t_{-})^{\frac{4k}{n-2k}} |dx|^2 = (\ln \frac{|x|}{r_{-}})^{\frac{4k}{n-2k}} |dx|^2,
\]
and the metric $v^{-2}|dx|^2$ is defined on $r_{-}< |x| < r_{+}$, and as $|x| \to r_{-}$, it has the
above degeneracy.

\item If $2k = n$, then $\sqrt[k]{D(\xi)} \sim \sqrt[k]{h}$ as $\xi \to \infty$.
Using this and  \eqref{+k},
we conclude that $t \to \pm \infty$ as $\xi \to \infty$. In fact,
 $t \sim \pm \frac {\xi}{\sqrt{1+\sqrt[k]{|h|}}}$. Taking the $-$ sign, for instance, we obtain
\[
g \sim |x|^{2(\sqrt{1+\sqrt[k]{|h|}}-1)} |dx|^2,
\]
so the metric $v^{-2}|dx|^2$ is defined on $0 < |x| < r_{+}$, and as $|x| \to 0$, it has the
above degeneracy.

\item If $2k > n$, then $\sqrt[k]{D(\xi)} \sim \sqrt[k]{h} e^{\frac{n-2k}{k} \xi}$ as $\xi \to \infty$. 
Using this and  \eqref{+k},
we conclude that $t \to \pm \infty$ as $\xi \to \infty$. In fact, for some costant $c$,
\[
 \xi \pm t = c - \frac{\sqrt[k]{|h|}}{2} \frac{k}{2k-n} e^{-\frac{2k-n}{k}\xi} + h.o.t.
\]
as $\xi  \to \infty$, 
from which we conclude that
\[
v^{-2} = e^{-2c}\{1 + \sqrt[k]{|h|} \frac{k}{2k-n} e^{-\frac{2k-n}{k}c}|x|^{\frac{2k-n}{k}}+ h.o.t. \}
\]
as $|x| \to 0$. Thus the metric $v^{-2}|dx|^2$ defined on $0 < |x| < r_{+}$ 
and extends to be a $C^{2-\frac{n}{k}}$ metric on $|x| < r_{+}$.
\end{enumerate}
\end{proof}
\begin{proof}[Proof of Theorem 3] When $\sigma_k(A_g) \equiv 0$, any radial solution of \eqref{3}
satisfies either $\xi_t^2 =1$ or $\xi_{tt} = (1-\frac{n}{2k})(1-\xi_t^2)$. Either case can be
integrated out easily. For instance, in the secend case, set $\eta = \xi_t$, then 
$\eta_t = (1-\frac{n}{2k})(1-\eta^2)$. It follows  that either $\eta \equiv \pm 1$ or 
\[
\eta = \frac {a e^{(2-\frac{n}{k})t}  +1}{a e^{(2-\frac{n}{k})t} -1}
\]
for some (non-zero) constant $a$. Integrating one more time concludes Theorem 3.
\end{proof}

\end{document}